\documentclass{elsart}%
\usepackage{amsfonts}
\usepackage{amsmath}
\usepackage{amssymb}
\usepackage{graphicx}%
\begin{document}

\date{22 September 2006}

%

%TCIMACRO{\TeXButton{Begin frontmatter}{\begin{frontmatter}}}%
%BeginExpansion
\begin{frontmatter}%
%EndExpansion
%

%TCIMACRO{\TeXButton{Title}{\title
%{Conditions for stability and instability of retrial queueing systems with general retrial times}%
%}}%
%BeginExpansion
\title
{Conditions for stability and instability of retrial queueing systems with general retrial times}%
%EndExpansion
%

%TCIMACRO{\TeXButton{Author}{\author{Tewfik Kernane}}}%
%BeginExpansion
\author{Tewfik Kernane}%
%EndExpansion
%

%TCIMACRO{\TeXButton{Address}{\address
%{Department of Probability and Statistics, Faculty of Mathematics,
%University of Sciences and Technology USTHB, Algiers, Algeria}}}%
%BeginExpansion
\address{Department of Probability and Statistics, Faculty of Mathematics,
University of Sciences and Technology USTHB, Algiers, Algeria}%
%EndExpansion
%

%TCIMACRO{\TeXButton{Address}{\address{E-mail: tkernane@gmail.com}}}%
%BeginExpansion
\address{E-mail: tkernane@gmail.com}%
%EndExpansion
%

%TCIMACRO{\TeXButton{Begin abstract}{\begin{abstract}} }%
%BeginExpansion
\begin{abstract}
%EndExpansion
We study the stability of single server retrial queues under general
distribution for retrial times and stationary ergodic service times, for three
main retrial policies studied in the literature: classical linear, constant
and control policies. The approach used is the renovating events approach to
obtain sufficient stability conditions by strong coupling convergence of the
process modeling the dynamics of the system to a unique stationary ergodic
regime. We also obtain instability conditions by convergence in distribution
to improper limiting sequences.
%TCIMACRO{\TeXButton{End abstract}{\end{abstract}}}%
%BeginExpansion
\end{abstract}%
%EndExpansion
%

%TCIMACRO{\TeXButton{Begin keywords}{\begin{keyword}} }%
%BeginExpansion
%EndExpansion
%

%TCIMACRO{\TeXButton{End frontmatter}{\end{frontmatter}}}%
%BeginExpansion
\end{frontmatter}%
%EndExpansion

\section*{Introduction}

The analysis of stability in queueing systems is the first step in studying
such models. The steady state solutions and performance characteristics of the
system do not exist if it is not stable. The efficiency of a queueing system
is related closely to its stability and is considered as inefficient if it is
unstable. Retrial queues have the characteristic that an arriving customer who
finds all waiting positions and service zones occupied must join a group of
"blocked" customers in an additional queue called "orbit" and reapplies for
getting served after random time intervals according to a specific retrial
policy. They arise in many practical situations. The classical example can be
found in telephone traffic theory where subscribers redial after receiving a
busy signal. For computer and communication applications, peripherals in
computer systems may make retrials to receive service from a central
processor. Another example can be adopted from the aviation where an aircraft
is directed into the waiting zone, if the runway is found busy, from which the
demand of landing is repeated at random periods of time. Retrial queueing
models are generally more complicated than traditional ones especially when
dealing with general distribution for retrial times. The existence of this
supplementary flow from the orbit and the random access to the server (as for
the linear policies that depend on the number of customers in orbit) make the
system more congested and difficult to model by simple random processes like
Markovian ones, which have properties that allow to derive easily conditions
for stability, especially when we do not assume an exponential distribution
(which has the memorylless property) allowing to obtain Markovian processes
modeling the system. Furthermore, it has been observed in telecommunication
systems that the exponential law is not a good estimator for the distribution
of retrial times (see Yang et al., 1994).

The subject of this paper is to analyze the stability of single server retrial
queues under general distribution for retrial times and stationary ergodic
service times (without independence assumption), for three main retrial
policies studied in the literature: classical linear, constant and control
policies. Stability results for such models with general retrial times are
rare and generally reduced to Markovian assumptions. For the linear retrial
policy, Koba and Kovalenko (2004) obtained a sufficient stability condition
(arrival rate is less than the service rate) for an M/G/1 system with
non-lattice distribution for retrial times satisfying an additional estimate
condition, with i.i.d service times. For the constant retrial policy, Koba
(2002) derived a stability condition for a GI/G/1 retrial system with a FIFO
discipline for the access from the orbit to the server and a general
distribution for orbit time in latticed and non-latticed cases with i.i.d
service times. For the control policy, Gomez-Corral (1999) studied extensively
an M/G/1 retrial queue with general retrial times where he derived the
stability condition for i.i.d service times and a FIFO discipline. For
non-independent service times, Altman and Borovkov (1997) obtained a
sufficient condition for the stability of a linear retrial queue under general
stationary ergodic service times and independent and exponentially distributed
interarrival and retrial times using the method of renovation events. Kernane
and A\"{\i}ssani (2006) obtained sufficient conditions for the stability of
various retrial queues with versatile retrial policy which incorporates the
constant and linear retrial policies under general stationary ergodic service
times and independent and exponentially distributed interarrival and retrial times.

The main approach used in this paper is the method of renovation events
originated in the work of Akhmarov and Leont'eva (1976) and developed by
Borovkov (1984) in the stationary ergodic setting. In the following section,
we derive stability and instability conditions for the classical linear
retrial policy with general retrial times, stationary ergodic service times
and Poisson arrivals. In Section 3, we obtain a stability condition and an
instability one for the constant retrial policy system with general retrial
times, stationary ergodic service times and Poisson arrivals. With the later
assumptions, we derive in Section 4, stability and instability conditions for
the control policy retrial model.

\section{Linear Retrial Policy}

We begin by considering the classical single server retrial system with linear
retrial policy. Customers arrive from outside according to a Poisson process
with rate $\lambda.$ If an arriving customer finds the server busy, he joins
the orbit and repeats his attempt to get served after random time intervals.
We consider the linear retrial policy where each customer in orbit attempts to
get served independently of other customers and we assume that the sequence of
inter-retrial times of a single customer is an independent sequence with
general distribution $R(\cdot)$, density function $r(\cdot)$ and Laplace
transform $r^{\ast}(z),$ $z>0.$ The successive service times $\left\{
\sigma_{n}\right\}  $ are assumed to form a stationary (in the strict sense)
and ergodic (which essentially means that time averages converge to constants
\textit{a.s}) sequence with $0<\mathbb{E}\sigma_{n}<\infty$. The
inter-arrival, inter-retrial and service times are assumed to be mutually independent.

Let $Q(t)$ be the number of customers in orbit at time $t$ and denote by
$s_{n}$ the instant when the $n$th service time ends. Consider the embedded
process $Q_{n}=Q(s_{n}+)$ of the number of customers in orbit just after the
end of the $n$th service duration. Denote by $N_{\lambda}(t)$ the counting
Poisson process with parameter $\lambda$ which counts the number of arriving
customers during a time interval $(0,t].$ If $Q_{n}=k,$ then we denote by
$\pi_{1}(n),...,\pi_{k}(n)$ the residual retrial times (forward recurrence
times) of the customers in orbit just after the instant $s_{n}$ and by
$\gamma_{n}$ the residual external arrival time at the same instant.

It is easy to see that the process $Q_{n}$ satisfies the following recurrence
relation:%
\begin{equation}
Q_{n+1}=(Q_{n}+\xi_{n})^{+}, \label{srs1}%
\end{equation}
where $x^{+}=\max[0,x]$ and%
\begin{equation}
\xi_{n}=N_{\lambda}(\sigma_{n})-\mathbb{I}\left\{  \min(\pi_{1}(n),...,\pi
_{Q_{n}}(n))<\gamma_{n}\right\}  , \label{main}%
\end{equation}
We have then expressed $Q_{n}$ as a Stochastic Recursive Sequence (SRS) (for
the definition see Borovkov, 1998).

We introduce the $\sigma-$algebra $\mathcal{F}_{n}^{\sigma}$ generated by the
set of random variables $\left\{  \sigma_{k}:k\leq n\right\}  $ and
$\mathcal{F}^{\sigma}$ generated by the entire sequence $\left\{  \sigma
_{n}:-\infty<n<+\infty\right\}  $ and for which any independent sequence not
depending on $\left\{  \sigma_{n}\right\}  $ is $\mathcal{F}^{\sigma}%
$-measurable (see Borovkov (1976) p.14). Let $U$ be the measure preserving
shift transformation of $\mathcal{F}^{\sigma}$-measurable random variables,
that is $U\sigma_{k}=\sigma_{k+1}$, and if $\eta\in\mathcal{F}^{\sigma}$ then
the sequence $\{\eta_{n}=U^{n}\eta:-\infty<n<+\infty\}$ is a stationary
ergodic sequence where $U^{n}$ is the $n$th iteration of $U$ and $U^{-n}$ is
the inverse transformation of $U^{n}$ $n\in%
%TCIMACRO{\U{2124} }%
%BeginExpansion
\mathbb{Z}
%EndExpansion
.$ We shall denote by $T$ the corresponding transformation of events in
$\mathcal{F}^{\sigma},$ that is for any $\mathcal{F}^{\sigma}$-measurable
sequence $\eta_{n}:$%
\begin{equation}
T\left\{  \omega:(\eta_{0}(\omega),...,\eta_{k}(\omega))\in(B_{0}%
,...,B_{k})\right\}  =\left\{  \omega:(\eta_{1}(\omega),...,\eta_{k+1}%
(\omega))\in(B_{0},...,B_{k})\right\}  ,
\end{equation}
where the events $B_{i}\in\mathcal{F}^{\sigma},$ $i=0,...,k.$

An event $A\in\mathcal{F}_{n+m}^{\xi},$ $m\geq0$, is a renovation event for
the SRS $\left\{  Q_{n}\right\}  $ on the segment $\left[  n,n+m\right]  $ if
there exists a measurable function $g$ such that on the set $A$%
\begin{equation}
Q_{n+m+1}=g(\xi_{n},...,\xi_{n+m}). \label{s1}%
\end{equation}
The sequence ${A_{n}}$, $A_{n}\in\mathcal{F}_{n+m}^{\xi}$, is a
\emph{renovating sequence of events} for the SRS $\left\{  Q_{n}\right\}  $ if
there exists an integer $n_{0}$ such that (\ref{s1}) holds true for $n\geq
n_{0}$ with a common function $g$ for all $n.$

We say that the SRS $\left\{  Q_{n}\right\}  $ is \textit{coupling convergent
}to a stationary sequence $\left\{  Q^{n}=U^{n}Q^{0}\right\}  $ if%
\begin{equation}
\lim_{n\rightarrow\infty}\mathbb{P}\left\{  Q_{k}=Q^{k};\text{ }\forall\text{
}k\geq n\right\}  =1. \label{cc}%
\end{equation}
Set $\nu_{k}=\min\{n\geq-k:$ $U^{-k}Q_{n+k}=Q^{n}\}$ and $\nu=\sup_{k\geq0}%
\nu_{k}.$

A SRS\ $\left\{  Q_{n}\right\}  $ is \textit{strong coupling convergent }to a
stationary sequence $\left\{  Q^{n}=U^{n}Q^{0}\right\}  $ if $\nu<\infty$ with
probability 1.

\begin{thm}
Assume that $\lambda\mathbb{E}\mathbf{\sigma}_{1}<1.$ Then the process
$\left\{  Q_{n}\right\}  $ is strong coupling convergent to a unique
stationary ergodic regime.\newline If $\lambda\mathbb{E}\mathbf{\sigma}%
_{1}>1,$ then the process $\left\{  Q_{n}\right\}  $ converges in distribution
to an improper limiting sequence.
\end{thm}

\begin{pf}
Since the driving sequence $\xi_{n}$ depend on $Q_{n},$ we will proceed first
by considering an auxiliary sequence $Q_{n}^{\ast}$ which majorizes $Q_{n}$
and having a driving sequence $\xi_{n}^{\ast}$ independent of $Q_{n}$ and it
has the following form:%
\begin{equation}
Q_{0}^{\ast}=Q_{0},\text{ \ \ }Q_{n+1}^{\ast}=\max(C,Q_{n}^{\ast}+\xi
_{n}^{\ast}),
\end{equation}
where%
\begin{equation}
\xi_{n}^{\ast}=N_{\lambda}(\sigma_{n})-\mathbb{I}\left\{  \min(\pi
_{1}(n),...,\pi_{C}(n))<\gamma_{n}\right\}  .
\end{equation}
The constant integer $C$ will be chosen later appropriately, and if
$Q_{n}^{\ast}>C$ the $C$ customers for which we consider the forward
recurrence times $\pi_{1}(n),...,\pi_{C}(n)$ are chosen randomly by an urn
scheme without repetition. Following the procedure used in Altman and Borovkov
(1997) and later in Kernane and A\"{\i}ssani (2006), we will construct
stationary renovation events with strictly positive probability for $Q_{n}$
from those of $Q_{n}^{\ast},$ and applying an ergodic theorem (Theorem 11.4 in
Borovkov, 1998) which states that an SRS is strong coupling convergent to a
unique stationary regime, satisfying the same recursion, if there exist
stationary renovating events of strictly positive probability.\newline The
stationarity and ergodicity of $\xi_{n}^{\ast}$ follows from the fact that
$\xi_{n}^{\ast}$ is $\mathcal{F}^{\sigma}$-measurable (for more details on the
ergodicity and stationarity of $\xi_{n}^{\ast}$ see Kernane and A\"{\i}ssani,
2006). We have%
\begin{align}
\mathbb{E}\xi_{n}^{\ast}  &  =\lambda\mathbb{E}\mathbf{\sigma}_{1}%
-\mathbb{P}\left(  \min(\pi_{1}(n),...,\pi_{C}(n))<\gamma_{n}\right) \\
&  =\lambda\mathbb{E}\mathbf{\sigma}_{1}-\left[  1-\mathbb{P}\left(  \pi
_{1}(n)\geq\gamma_{n},...,\pi_{C}(n)\geq\gamma_{n}\right)  \right] \\
&  =\lambda\mathbb{E}\mathbf{\sigma}_{1}-\left[  1-\int\limits_{0}^{\infty
}\lambda e^{-\lambda t}\text{ }\mathbb{P}\left(  \pi_{1}(n)\geq t,...,\pi
_{C}(n)\geq t\right)  \text{ }dt\right] \\
&  =\lambda\mathbb{E}\mathbf{\sigma}_{1}-\left[  1-\int\limits_{0}^{\infty
}\lambda e^{-\lambda t}\prod\limits_{i=1}^{C}\mathbb{P}(\pi_{i}(n)\geq
t)\text{ }dt\right]  .
\end{align}
Since $\lim_{C\rightarrow\infty}\prod\limits_{i=1}^{C}\mathbb{P}(\pi
_{i}(n)\geq t)=0,$ then by dominated convergence theorem%
\begin{equation}
\lim_{C\rightarrow\infty}\int\limits_{0}^{\infty}\lambda e^{-\lambda t}%
\prod\limits_{i=1}^{C}\mathbb{P}(\pi_{i}(n)\geq t)\text{ }dt=0.
\end{equation}
If the condition $\lambda\mathbb{E}\mathbf{\sigma}_{1}<1$ is satisfied, we can
choose the constant $C$ such that $\mathbb{E}\xi_{n}^{\ast}<0.$ It follows
from example 11.1 in Borovkov (1998) that there exists a stationary renovating
sequence of events with positive probability for $Q_{n}^{\ast},$ from which we
deduce those of $Q_{n}$ (see Altman and Borovkov, 1997). Applying the ergodic
theorem (Theorem 11.4 in Borovkov, 1998) we obtain that the sequence $Q_{n}$
is strong coupling convergent to a unique stationary process $\widetilde
{Q}_{n}=U^{n}\widetilde{Q}_{0},$ with $\widetilde{Q}_{0}$ $\mathcal{F}%
^{\sigma}$-measurable and since $\widetilde{Q}_{n}$ is an $U-$shifted
$\mathcal{F}^{\sigma}$-measurable sequence then it is ergodic.\newline For the
instability condition, consider the auxiliary process $Q_{n}^{S}$
corresponding to a simple single server queue without retrials, that is%
\begin{equation}
Q_{0}^{S}=Q_{0},\text{ \ \ }Q_{n+1}^{S}=(Q_{n}^{S}+\xi_{n}^{S})^{+},
\end{equation}
where%
\begin{equation}
\xi_{n}^{S}=N_{\lambda}(\sigma_{n})-1.
\end{equation}
Clearly $Q_{n}^{S}\leq_{st}Q_{n}$ and it is well known that if $\lambda
\mathbb{E}\mathbf{\sigma}_{1}>1$ then $\lim_{n\rightarrow\infty}Q_{n}%
^{S}=+\infty$ \textit{a.s. }(see Theorem 1.7 in Borovkov, 1976).\textit{
}Thus, the process $\left\{  Q_{n}\right\}  $ converges in distribution to an
improper limiting sequence.
\end{pf}

\section{Constant Retrial Policy}

Consider now a single server retrial queue governed by the constant retrial
policy which is described as follows. After a random time generally
distributed (which we will call the orbit retrial time), one customer from the
orbit (at the head of the queue or a randomly chosen one if any) take his
service if the server is free, so an orbit time can be in progress even though
the server is busy, this may happen in system where the orbit has no
information about the state of the server. The sequence of orbit cycle times
$\left\{  r_{i}\right\}  $ is assumed to be i.i.d, having $R(\cdot)$ as cdf,
$r(\cdot)$ as density function with mean $\mathbb{E}r_{1}$ and Laplace
transform $r^{\ast}(z),$ $z>0.$ Let $\pi(n)$ be the forward recurrence time of
the orbit retrial time after the end of the $n$th service time. Then the
process $Q_{n}$ has now the following representation as a SRS:%
\begin{equation}
Q_{n+1}=(Q_{n}+\xi_{n})^{+},
\end{equation}
where%
\begin{equation}
\xi_{n}=N_{\lambda}(\sigma_{n})-\mathbb{I}\left\{  \pi(n)<\gamma_{n}\right\}
.
\end{equation}

\begin{thm}
If $R$ is nonlattice and%
\begin{equation}
\lambda\mathbb{E}\mathbf{\sigma}_{1}<\dfrac{\left[  1-r^{\ast}(\lambda
)\right]  }{\lambda\mathbb{E}r_{1}}, \label{const-cond}%
\end{equation}
then the process $\left\{  Q_{n}\right\}  $ is strong coupling convergent to a
unique stationary ergodic regime.\newline If $\lambda\mathbb{E}\sigma
_{1}>(1-r^{\ast}(\lambda))/(\lambda\mathbb{E}r_{1})$. Then the process $Q_{n}$
converges in distribution to an improper limiting sequence.
\end{thm}

\begin{pf}
We have%
\begin{equation}
\mathbb{E}\xi_{n}=\lambda\mathbb{E}\mathbf{\sigma}_{1}-\mathbb{P}\left(
\pi(n)<\gamma_{n}\right)  .
\end{equation}
Since the interarrival times are exponentially distributed then so is the
residual arrival time $\gamma_{n},$ hence%
\begin{equation}
\mathbb{P}\left(  \pi(n)<\gamma_{n}\right)  =\int\limits_{0}^{+\infty
}\mathbb{P}\left(  \pi(n)<t\right)  \text{ }\lambda\text{ }e^{-\lambda t}dt.
\label{pin-gam}%
\end{equation}
Since we are interesting on steady state behaviour of the system and by
assuming a nonlattice (also called non-arithmetic) distribution $R(t)$ for
orbit retrial times, then from a well known result in renewal theory (see Cox,
1962) we have the following asymptotic distribution for the forward recurrence
time
\begin{equation}
\mathbb{P}\left(  \pi(n)<t\right)  =\frac{1}{\mathbb{E}r_{1}}\int
\limits_{0}^{t}\left[  1-R(x)\right]  \text{ }dx.
\end{equation}
The formula (\ref{pin-gam}) becomes%
\begin{align}
\mathbb{P}\left(  \pi(n)<\gamma_{n}\right)   &  =\frac{1}{\mathbb{E}r_{1}}%
\int\limits_{0}^{+\infty}\left[  1-R(x)\right]  \text{ }\int\limits_{x}%
^{+\infty}\lambda e^{-\lambda t}dt\text{ }dx\\
&  =\frac{1}{\mathbb{E}r_{1}}\int\limits_{0}^{+\infty}\left[  1-R(x)\right]
\text{ }e^{-\lambda x}\text{ }dx=\frac{1}{\mathbb{E}r_{1}}\left[
\frac{1-r^{\ast}(\lambda)}{\lambda}\right]  .
\end{align}
Now if condition (\ref{const-cond}) is satisfied then $\mathbb{E}\xi_{n}<0.$
Since $\xi_{n}$ is $\mathcal{F}^{\sigma}$-measurable (generated by $\sigma
_{n}$) then it is a stationary ergodic sequence. From this and example 11.1 in
Borovkov (1998), there exists a stationary sequence of renovation events with
positive probability for $\{Q_{n}\}$. Hence, using Theorem 11.4 of Borovkov
(1998), the sequence $\{Q_{n}\}$ is strong coupling convergent to a unique
stationary sequence $\widetilde{Q}_{n}$ obeying the equation $\widetilde
{Q}_{n+1}=(\widetilde{Q}_{n}+\xi_{n})^{+},$ the ergodicity of $\widetilde
{Q}_{n}$ follows from the fact that $\widetilde{Q}_{n}$ is an $U-$shifted
sequence ($\widetilde{Q}_{n}=U^{n}\widetilde{Q}_{0},$ with $\widetilde{Q}_{0}$
$\mathcal{F}_{0}^{\sigma}$-measurable) generated by the stationary and ergodic
sequence $\xi_{n}$.\newline The instability condition $\lambda\mathbb{E}%
\sigma_{1}>\left[  1-r^{\ast}(\lambda)\right]  /\lambda\mathbb{E}r_{1}$ yields
to $\mathbb{E}\xi_{n}>0,$ and it is well known that for SRS of the form
$Q_{n+1}=(Q_{n}+\xi_{n})^{+}$ this implies the convergence of the process
$Q_{n}$ to an improper limiting sequence (see Theorem 1.7 of Borovkov (1976)).
\end{pf}

\subsection{Exponential retrial times}

By assuming an exponential distribution with parameter $\theta$ for retrial
times, that is $R(x)=1-e^{-\theta x},$ it is well known that $r^{\ast
}(s)=s/(s+\theta),$ and $\mathbb{E}r_{1}=1/\theta.$ The condition
(\ref{const-cond}) will read up, after some algebra, as follows%
\begin{equation}
\lambda\mathbb{E}\mathbf{\sigma}_{1}<\dfrac{\theta}{\lambda+\theta}.
\end{equation}
Which is the condition obtained in the paper of Kernane and A\"{\i}ssani
(2006), in exponential retrial context.

\section{Retrial Control Policy}

Consider a single server retrial queue with a control retrial policy. Primary
customers enter from the outside according to a Poisson process with rate
$\lambda.$ If a primary customer finds the server busy upon arrival it joins
the orbit to connect later according to the control retrial policy, which is
described as follows. Just after the end of a service time a generally
distributed retrial time begins to find the server free. If the retrial time
finishes before an external arrival, then one customer from the orbit (at the
head of the queue or a randomly chosen one if any) receives its service and
leaves the system. We assume that the sequence of retrial times $\left\{
r_{n}\right\}  $ is an i.i.d sequence having $r(\cdot)$ as pdf$,$ $R(\cdot)$
as cdf and Laplace transform $r^{\ast}(\cdot)$, with finite mean
$\mathbb{E}r_{1}$. The $n$th service duration of a call is $\sigma_{n}$ and we
assume that the sequence of service times $\left\{  \sigma_{n}\right\}  $ is
stationary and ergodic with $0<\mathbb{E}\sigma_{1}<\infty$.

The process $\left\{  Q_{n}\right\}  $ has the following representation as a
stochastic recursive sequence SRS:%
\begin{equation}
Q_{n+1}=(Q_{n}+\xi_{n})^{+},
\end{equation}
where%
\begin{equation}
\xi_{n}=N_{\lambda}(\sigma_{n})-\mathbb{I}\left\{  r_{n}<\gamma_{n}\right\}  ,
\end{equation}
where $\gamma_{n}$ is the residual arrival time of an external call at the end
of the $n$th service period.

\begin{thm}
Assume that
\begin{equation}
\lambda\mathbb{E}\sigma_{1}<r^{\ast}(\lambda). \label{cont-policy}%
\end{equation}
Then the process $Q_{n}$ is strong coupling convergent to a unique stationary
ergodic regime.\newline If $\lambda\mathbb{E}\sigma_{1}>r^{\ast}(\lambda)$.
Then the process $Q_{n}$ converges in distribution to an improper limiting sequence.
\end{thm}

\begin{pf}
The proof is similar to that of Theorem 2, by noting that $\mathbb{E}\xi
_{n}=\lambda\mathbb{E}\sigma_{1}-r^{\ast}(\lambda).$
\end{pf}

\subsection{Exponential retrial times}

Assume that the retrial times are exponentially distributed with mean
$1/\theta,$ then $r^{\ast}(\lambda)=\theta/(\lambda+\theta)$ and the stability
condition (\ref{cont-policy}) becomes:%
\begin{equation}
\lambda\mathbb{E}\sigma_{1}<\frac{\theta}{\lambda+\theta}.
\end{equation}
This condition is quite evident since it can be obtained from the constant
policy from the memorylless property of the exponential distribution.

\subsection{Hyperexponential distribution for retrial times}

Assume now that the retrial times follow the hyperexponential distribution
with density $r(x)=p\theta\exp(-\theta x)+\left(  1-p\right)  \theta^{2}%
\exp(-\theta^{2}x),$ $0\leq p<1.$ Then $r^{\ast}(\lambda)=\theta\left[
\lambda\left(  p+\left(  1-p\right)  \theta\right)  +\theta^{2}\right]
/(\lambda+\theta)\left(  \lambda+\theta^{2}\right)  ,$ and the stability
condition (\ref{cont-policy}) in this case is%
\begin{equation}
\lambda\mathbb{E}\sigma_{1}<\frac{\theta\left[  \lambda\left(  p+\left(
1-p\right)  \theta\right)  +\theta^{2}\right]  }{(\lambda+\theta)\left(
\lambda+\theta^{2}\right)  }.
\end{equation}

\subsection{The Erlang distribution for retrial times}

The Erlang distribution has been found useful for describing random variables
in queueing applications. The density of an $Erlang\left(  n,\mu\right)  $
distribution is given by $r(x)=\mu^{n}\exp(-\mu x)x^{n-1}/\left(  n-1\right)
!,$ $x>0$ and $n\in%
%TCIMACRO{\U{2115} }%
%BeginExpansion
\mathbb{N}
%EndExpansion
^{\ast}.$ Its Laplace transform is $r^{\ast}(s)=\mu^{n}/\left(  s+\mu\right)
^{n}.$ Then the control policy model will be stable if%
\begin{equation}
\lambda\mathbb{E}\sigma_{1}<\left(  \frac{\mu}{\lambda+\mu}\right)  ^{n}.
\end{equation}

\begin{rem}
It should be noted that the assumption $\mathbb{E}\xi_{n}=0,$ and weak
dependence among the $\xi_{n},$ does not preclude the possibility that the
process $\{Q_{n}\}$ converges to a proper stationary regime.
\end{rem}

\begin{rem}
Conditions for the stability of modified models with general retrial times,
such as allowing breakdowns of the server, two types of arrivals, negative
arrivals and batch arrivals models may be obtained easily following the
procedure used in Kernane and A\"{\i}ssani (2006). The conditions of the
stability will be written by replacing the left hand side of the classical
linear policy by the left hand sides of the case of a linear versatile policy
obtained in Kernane and A\"{\i}ssani (2006). For the constant policy, we have
to make the appropriate changes to the driving sequences in the SRS modeling
the dynamics of the modified models in Kernane and A\"{\i}ssani (2006) by
considering the residual orbit time as shown here in Section 3, the conditions
of stability will follow directly after some algebra.
\end{rem}

\begin{rem}
We may also consider the versatile retrial policy by incorporating the
residual orbit retrial time $\pi(n)$ in the equation \ref{main} and
considering the whole retrial times of the customers in orbit $r_{1}%
(n),...,r_{Q_{n}}(n)$ as follows%
\begin{equation}
\xi_{n}=N_{\lambda}(\sigma_{n})-\mathbb{I}\left\{  \min(\pi(n)+r_{1}%
(n),...,\pi(n)+r_{Q_{n}}(n))<\gamma_{n}\right\}  ,
\end{equation}
the condition of Theorem 1 still holds for this versatile retrial policy, by
noting that in the proof we have to consider
\begin{equation}
\mathbb{E}\xi_{n}^{\ast}=\lambda\mathbb{E}\mathbf{\sigma}_{1}-\left[
1-\int\limits_{0}^{\infty}\int\limits_{0}^{t}\lambda e^{-\lambda t}\text{
}P\left(  r_{1}(n)\geq t-s,...,r_{C}(n)\geq t-s\right)  \text{ }%
dG(s)dt\right]  ,
\end{equation}
where $G(s)$ is the cdf of the residual orbit retrial time satisfying%
\begin{equation}
G(s)=\frac{1}{\mathbb{E}\alpha_{1}}\int\limits_{0}^{s}\left[  1-A(x)\right]
\text{ }dx,
\end{equation}
with $\mathbb{E}\alpha_{1}$ the mean of the orbit retrial time and $A(x)$ its cdf.
\end{rem}

\end{document}